\documentclass{mcom-l}

\usepackage{amssymb}		
\usepackage{eucal}		
\usepackage{lscape}		
\usepackage{longtable}		

\renewcommand\ge\geqslant
\renewcommand\le\leqslant

\newtheorem{theorem}{Theorem}[section]
\newtheorem{corollary}[theorem]{Corollary}
\newtheorem{proposition}[theorem]{Proposition}

\theoremstyle{definition}
\newtheorem{definition}[theorem]{Definition}
\newtheorem{example}[theorem]{Example}

\theoremstyle{remark}

\newtheorem{remarkaftertheorem}{Remark}[theorem]

\numberwithin{equation}{section}

\newcommand{\defeq}{\stackrel{\rm{def}}{=}}
\newcommand{\Hl}{\mathop{{}\it Hl}\nolimits}

\begin{document}

\title{The 192 Solutions of the Heun Equation}

\author{Robert S. Maier}
\address{Depts.\ of Mathematics and Physics, University of Arizona, Tucson,
AZ 85721, USA}
\email{rsm@math.arizona.edu}
\urladdr{http://www.math.arizona.edu/\~{}rsm}
\thanks{The author was supported in part by NSF Grant No.\ PHY-0099484.}

\subjclass[2000]{Primary 33E30; 33-04, 34M15, 33C05, 20F55.}
\date{}

\begin{abstract}
A machine-generated list of $192$ local solutions of the Heun equation is
given.  They are analogous to Kummer's $24$~solutions of the Gauss
hypergeometric equation, since the two equations are canonical Fuchsian
differential equations on the Riemann sphere with four and three singular
points, respectively.  Tabulation is facilitated by the identification of
the automorphism group of the equation with $n$~singular points as the
Coxeter group~$\mathcal{D}_n$.  Each of the $192$ expressions is labeled by
an element of~$\mathcal{D}_4$.  Of~the $192$, $24$~are equivalent
expressions for the local Heun function~$\Hl$, and it~is shown that the
resulting order-$24$ group of transformations of~$\Hl$ is isomorphic to the
symmetric group~$S_4$.  The isomorphism encodes each transformation as a
permutation of an abstract four-element set, not identical to the set of
singular points.
\end{abstract}

\maketitle

\section{Introduction}

The theory of second-order linear differential equations on the Riemann
sphere $\mathbb{P}^1(\mathbb{C})\ni x$, with coefficients in the field of
rational functions~$\mathbb{C}(x)$, has been developed in both conceptual
and algorithmic directions.  In~a modern interpretation, any such equation
specifies a flat meromorphic connection on a two-dimensional vector bundle
over~$\mathbb{P}^1(\mathbb{C})$~\cite{Varadarajan96}.  The study of these
equations and their solutions led to the classical theory of special
functions~\cite{Poole36}.

The simplest, but still nontrivial case is the Fuchsian one, in which all
singular points of the equation are regular.  If~there are $n$~regular
singular points on~$\mathbb{P}^1(\mathbb{C})$ ($n\ge3$), one of which is
allowed to be~$\infty$, the equation will be determined up~to isomorphism
by their locations, a~parameter associated with each singular point (i.e.,
a~difference of characteristic exponents), and $n-3$~accessory parameters,
which have more a global than a local significance.  The case $n=3$ is the
simplest.  By appropriate normalizations, any Fuchsian equation with three
singular points can be reduced to the Gauss hypergeometric equation
\begin{equation}
\label{eq:hyperg}
D_x^2\, u + \left[\frac{c}{x} + \frac{a+b-c+1}{x-1}\right]
D_x u + \left[\frac{ab}{x(x-1)}\right]u = 0,
\end{equation}
which has singular points $x=0,1,\infty$, and respective characteristic
exponents ${0,1-c};\allowbreak\,{0,-a-b+c};\allowbreak\,{a,b}$.  The
analysis of~(\ref{eq:hyperg}) led to much interesting nineteenth-century
mathematics~\cite{Gray2000}.  Kummer's $24$~local solutions of this
equation, reproduced in Table~\ref{tab:hyperg} below, are especially well
known
\cite{Abramowitz64,Dwork84,Erdelyi53,Gray2000,Poole36,Prosser94,Whittaker27}.
Each is expressed in~terms of the Gauss hypergeometric function~${}_2F_1$.
The $24$~solutions split into $6$~sets of $4$~formally distinct but
equivalent expressions, each set defining one of the two Frobenius
solutions in the neighborhood of a singular point.

The classical study of the hypergeometric equation was facilitated by its
having no~accessory parameters, and the situation when $n>3$ is
significantly more complicated.  It is known that any Fuchsian equation
with $n=4$ singular points can be reduced to the Heun
equation~\cite{Ronveaux95}.  This is the equation
\begin{equation}
\label{eq:Heun}
D_x^2\, u
+ \left[ \frac\gamma x + \frac\delta{x-1} + \frac\epsilon{x-a}
  \right]D_x u + \left[\frac{\alpha\beta\, x - q}{x(x-1)(x-a)}\right]u = 0,
\end{equation}
which has singular points $x=0,1,a,\infty$, where
$a\in\mathbb{C}\setminus\{0,1\}$~is a free parameter.  The respective
exponents are ${0,1-\gamma};\allowbreak\,{0,1-\delta};\allowbreak\,
{0,1-\epsilon};\allowbreak\,\alpha,\beta$.  It follows from Fuchs's
relation on the exponents of Fuchsian equations that
$\epsilon=\alpha+\beta-\gamma-\delta+1$.  (If~this condition does not hold,
the singular point $x=\infty$ will not be regular.)  The Heun equation
includes an accessory parameter, namely the quantity $q\in\mathbb{C}$,
which in many applications appears as a spectral parameter.  This equation
is under active study~\cite{Maier03}, in~part because the band structure of
its solutions is important to the theory of integrable nonlinear wave
equations~\cite{Gesztesy95a,Smirnov2001}.  Fuchsian equations with $n=5$
singular points are also being systematically studied~\cite{Smirnov2003}.

Kummer's family of $24$ solutions of~(\ref{eq:hyperg}) has an analogue for
each $n\ge3$, namely a family of $2^{n-1}n!$ local solutions, which splits
into $2n$ sets of $2^{n-2}(n-1)!$~equivalent expressions, each set defining
one of the two Frobenius solutions in the neighborhood of a singular point.
(The $n!$ factor comes from permuting the $n$~singular points, and the
$2^{n-1}$ factor from negating exponent differences.)  It~is quite
difficult to write~down the full family of $2^{n-1}n!$ solutions explicitly
when~$n>3$, especially if one uses hand computation.  The $n=4$ case is a
case in~point.  Just as each of Kummer's $24$~solutions can be written
in~terms of the Gauss function~${}_2F_1$, each of the
$2^3\cdot4!=192$~solutions of the Heun equation can be written in~terms of
a canonical `local Heun function'~$\Hl$.  But no~satisfactory list of these
solutions has been published, until~now.  The stumbling block is the
calculation, in~each, of the value of the accessory parameter of~$\Hl$.
Heun listed $48$ of the~$192$ in his 1889 paper~\cite{Heun1889}, with the
accessory parameter omitted in each.  Applied mathematicians have noticed
that apart from this, his list contains misprints, and cannot be used in
practical applications~\cite[\S\,6.3]{Schmitz94}.  As an alternative the
work of Snow~\cite{Snow52} is frequently cited, since it gives $25$~of
the~$192$, in a slightly individualistic notation.  Several solutions are
given in the useful review of the late
F.~Arscott~\cite[Part~A]{Ronveaux95}, which also discusses how they can be
transformed into one another.

In this paper we determine the {\em group structure\/} of the set of
transformations that can be applied to any normalized Fuchsian equation
on~$\mathbb{P}^1(\mathbb{C})$ with $n$~singular points, to generate
alternative expressions for its family of local solutions.  This
automorphism group has order $2^{n-1}n!$ and acts on the parameter space of
the equation.  As Theorem~\ref{thm:asymmetric} states, it is isomorphic to
the Coxeter group~$\mathcal{D}_n$~\cite{Grove85}, the group of even-signed
permutations of an $n$-set.  This is an index-$2$ normal subgroup of the
hyperoctahedral Coxeter group~$\mathcal{B}_n$, which is the group of signed
permutations of an $n$-set.  ($\mathcal{B}_n$~itself is isomorphic to the
wreath product $\mathbb{Z}_2 \wr S_n$, where $S_n$ is the symmetric group.)
A~similar result was recently obtained by Oblezin~\cite{Oblezin2004}, in an
abstract framework.  That the automorphism group
is~$\mathcal{D}_n\!\vartriangleleft\mathcal{B}_n=[\mathcal{D}_n\!:\!\mathbb{Z}_2]$
clarifies an observation of Dwork~\cite{Dwork84} that the `Kummer
automorphism group', the order-$24$ automorphism group of the Gauss
hypergeometric equation~($n=3$), is isomorphic to the octahedral rotation
group~$S_4$, and has a natural $\mathbb{Z}_2$-extension to an order-$48$
group.  It~is an accident that this order-$24$ group is isomorphic to a
symmetric group: there is no~counterpart to the isomorphism $\mathcal{D}_3
\cong S_4$ when~$n>3$.

The identification of the automorphism group of the Heun equation
as~$\mathcal{D}_4$ facilitates both hand and machine computation of the
$192$~solutions, including accessory parameter values.  Only a few
generating transformations need to be repeatedly applied.  The solutions,
indexed by $g\in\mathcal{D}_4$, are presented in Table~\ref{tab:heun}.
It~unfortunately follows from this table that nearly all the 48~solutions
generated by Heun~\cite{Heun1889} contain at~least one error.  The errors
are more serious than mere misprints.

Interestingly, the $24$~solutions equivalent to~$\Hl$ (including
$\Hl$~itself) yield transformations of~$\Hl$ which under composition, form
an order-$24$ subgroup of~$\mathcal{D}_4$ isomorphic to the Kummer
automorphism group~$\mathcal{D}_3$, and hence to~$S_4$.  The isomorphism
to~$S_4$ is given in Table~\ref{tab:iso}.  The possibility of encoding the
$24$~basic transformations of~$\Hl$ as the permutations of an abstract
four-element set, not to be identified with the set of four singular
points, seems not to have been noticed before, though it is useful in hand
computation.


The question of the global structure of the space of Fuchsian equations
with $n$~singular points, up~to isomorphism, is not addressed here.  The
space of three distinct points in~$\mathbb{P}^1(\mathbb{C})$ is effectively
a single point, since any two such configurations are related by a unique
M\"obius transformation; and the space of four distinct points
in~$\mathbb{P}^1(\mathbb{C})$ is effectively~$\mathbb{C}$, which
compactifies to~$\mathbb{P}^1({\mathbb{C}})$.  The global structure of the
compactified space of five-point configurations has only recently been
determined~\cite{Yoshida94}.  Supplementing the singular points by exponent
differences and accessory parameters makes global issues significantly more
difficult to treat.

\section{Basic Facts}
\label{sec:basic}

Any second-order Fuchsian equation on~$\mathbb{P}^1(\mathbb{C})$ with
$n-1$~distinct finite singular points $d_1,\dots,d_{n-1}\in\mathbb{C}$ can
be written as~\cite{Poole36}
\begin{equation}
\label{eq:fuchsian}
T u \defeq
D_x^2\,u + \left[\sum_{i=1}^{n-1} \frac{1-\rho_i-\hat\rho_i}{x-d_i} \right]D_xu
+ \left[
\sum_{i=1}^{n-1}
\frac{\rho_i\hat\rho_i}{(x-d_i)^2}
+
\frac{\Pi_{n-3}(x)}{\prod_{i=1}^{n-1}(x-d_i)}
\right]u = 0.
\end{equation}
Here $\rho_i,\hat\rho_i\in\mathbb{C}$ are the characteristic exponents at
the point~$d_i$, which can be calculated by the method of Frobenius as the
roots of an indicial equation.  $\Pi_{n-3}$~is a polynomial of
degree~$n-3$, the leading coefficient of which enters into the indicial
equation at the point~$x=\infty$.  The $2n$~exponents, including the two
exponents $\rho^{(\infty)}\!,{\hat\rho}^{(\infty)}$ at infinity, obey
Fuchs's relation: their sum is~$n-2$.  The $n-3$ trailing coefficients
of~$\Pi_{n-3}$ are independent of the exponents.  According to one
convention, these are the $n-3$ accessory parameters.

The automorphisms of~(\ref{eq:fuchsian}) include the group of M\"obius
transformations of the independent variable~$x$, i.e., the projective
linear group~$PGL(2,\mathbb{C})$.  If $x=Py\defeq(Ay+B)/(Cy+D)$ with
$AD-BC\neq0$, the equation $Tu=0$ for $u=u(x)$ can be pulled back to an
equation $\widetilde T\tilde u=0$ for $\tilde u=\tilde u(y)\defeq u(Py)$.
This has singular points $y=P^{-1}d_1,\dots,P^{-1}d_{n-1}$ (and also
$P^{-1}\infty$, if $\infty$~is a singular point of~$Tu=0$).  The exponents
are preserved under M\"obius transformations, though the accessory
parameters are not.  The exponents of $\widetilde T\tilde u=0$ at any
$y\in\mathbb{P}^1(\mathbb{C})$ will be the exponents of $Tu=0$ at~$x=Py$.
This applies to ordinary points as~well as singular points.  The exponents
of any finite ordinary point are~$0,1$, and if the point at infinity is
ordinary, its exponents will be~$0,-1$.

The automorphisms of~(\ref{eq:fuchsian}) also include a group of
transformations, linear in the dependent variable, called F-homotopic or
index transformations.  For each $k\ge1$, such a transformation is
specified by distinct finite points $e_1,\dots,e_k\in\mathbb{C}$ and
corresponding `index shifts' $s_1,\dots,s_k\in\mathbb{C}$.  The transformed
equation is $\widehat T\hat u=0$, where $\hat u=\hat u(x)\defeq S(x)u(x)$
and $\widehat T\defeq STS^{-1}$, with $S=S(x)\defeq\prod_{i=1}^k
(x-e_i)^{s_i}$.  By examination, $\widehat T$~has the same general form
as~$T$.  In $\widehat T\hat u=0$, the exponents of each point $x=e_i$ are
shifted by~$s_i$ relative to those in~$Tu=0$, and the exponents
of~$x=\infty$ are shifted by~$-\sum_{i=1}^ks_i$.  In~general, the accessory
parameters are also altered.  For each choice of $e_1,\dots,e_k$, the index
transformations form an abelian group isomorphic to the additive
group~$\mathbb{C}^k$.

The structure of the space of Fuchsian equations up~to index
transformations is illuminated by two normalized versions
of~(\ref{eq:fuchsian}), traditionally called reduced
forms~\cite[\S\,20]{Poole36}.  The first is the {\em symmetrically reduced
form\/}
\begin{equation}
\label{eq:fuchsian0}
T_0 u \defeq
D_x^2\,u 
+ \left[
\sum_{i=1}^{n-1}
\frac{\frac14(1-\delta_i^2)}{(x-d_i)^2}
+
\frac{\Pi^{0}_{n-3}(x)}{\prod_{i=1}^{n-1}(x-d_i)}
\right]u = 0,
\end{equation}
which has $\rho_i+\hat\rho_i=1$ for each~$i$ and
$\rho^{(\infty)}\!+{\hat\rho}^{(\infty)}=-1$.  It is obtained
from~(\ref{eq:fuchsian}) by the index transformation which for each~$i$,
shifts the pair $\rho_i,\hat\rho_i$ by $(1-\rho_i-\hat\rho_i)/2$, and
$\rho^{(\infty)}\!,{\hat\rho}^{(\infty)}$ by
$-\sum_{i=1}^{n-1}(1-\rho_i-\hat\rho_i)/2$.  In~(\ref{eq:fuchsian0}),
$\delta_i^2\defeq (\rho_i-\hat\rho_i)^2$.  The symmetrically reduced form
could be denoted $\mathcal{E}_{\rm
sym}(d_1,\delta_1;\dots;d_{n-1},\delta_{n-1};\infty,\delta^{(\infty)};Q)$,
where $Q$~is the vector of $n-3$ accessory parameters.

The other is the {\em asymmetrically reduced form\/}
\begin{equation}
\label{eq:fuchsian1}
T_1 u \defeq
D_x^2\,u + \left[\sum_{i=1}^{n-1} \frac{1-\rho_i}{x-d_i} \right]D_xu
+ \left[
\frac{\Pi^1_{n-3}(x)}{\prod_{i=1}^{n-1}(x-d_i)}
\right]u = 0,
\end{equation}
which has one exponent (by convention~$\hat\rho_i$) equal to zero at each
finite singular point.  It~is obtained from~(\ref{eq:fuchsian}) by an index
transformation which shifts each of the $n-1$ pairs $\rho_i,\hat\rho_i$ by
either $-\rho_i$ or~$-\hat\rho_i$; and also shifts
$\rho^{(\infty)}\!,{\hat\rho}^{(\infty)}$ correspondingly.  The reduction
is not unique, since in~general there are $2^{n-1}$~such index
transformations.  It~can be shown that the leading coefficient
of~$\Pi_{n-3}^1$, i.e., the coefficient of~$x^{n-3}$,
equals~$\rho^{(\infty)}{\hat\rho}^{(\infty)}$.  This reduced form could be
denoted $\mathcal{E}_{\rm
asym}(d_1,\rho_1;\dots;d_{n-1},\rho_{n-1};\infty,\rho^{(\infty)}\!,{\hat\rho}^{(\infty)};Q)$,
where only $n$ of the exponent parameters are independent, since they must
satisfy Fuchs's relation.  

Due~to the nonuniqueness of the reduction, an order-$2^{n-1}$ \!group of
index transformations acts on this reduced form.  It~is generated by the
commuting involutions
\begin{equation}
\label{eq:tired}
(d_i,\rho_i;\infty,\rho^{(\infty)}\!,{\hat\rho}^{(\infty)}) \mapsto
(d_i,-\rho_i;\infty,\rho_i+{\hat\rho}^{(\infty)}\!,\rho_i+\rho^{(\infty)}),
\end{equation}
$i=1,\dots,n-1$, and is isomorphic to~$(\mathbb{Z}_2)^{n-1}$.  On the level
of the unreduced Fuchsian equation~(\ref{eq:fuchsian}), these involutions
negate the exponent {\em differences\/} at the points
$x=d_1,\dots,d_{n-1}$.  It~should be noted that the exponents
$\rho^{(\infty)}\!,{\hat\rho}^{(\infty)}$ are determined
by~(\ref{eq:fuchsian1}) only up~to interchange; this seeming ambiguity in
the action~(\ref{eq:tired}) will be discussed further below.

How to quotient~out M\"obius transformations, in~addition to index
transformations, will also be discussed.  The concept of a {\em signed
permutation\/} of the $n$~singular points, where the signs keep track of
the negatings of exponent differences, will prove useful.  A~signed
permutation of a set of objects $O_n=\{1,\dots,n\}$ is a bijection
$\pi:O_n\to O_n$ in which each $i\mapsto\pi(i)$ is accompanied by a `sign'.
Equivalently, it~is a bijection $\pi':O_n'\to O_n'$, where
$O_n'=\{1,\dots,n\}\cup\{-1,\dots,-n\}$, that satisfies $\pi'(-i)=-\pi'(i)$
for~all~$i$.  The group of signed permutations of~$O_n$ has order~$2^n n!$
and is isomorphic to the wreath product~$\mathbb{Z}_2 \wr S_n$.  Any signed
permutation may be written in a sign-annotated version of the disjoint
cycle representation used for elements of~$S_n$.  For example,
$[1_+2_-3_+][4_-]$ signifies $1\mapsto2$, $2\mapsto-3$, $3\mapsto1$,
$4\mapsto-4$.  Positively signed $1$-cycles may be omitted.  A signed
permutation is said to be even-signed if its annotated cycle representation
has an even number of negative signs.  Even-signed permutations form a
normal subgroup of index~$2$, which is isomorphic to a semi-direct product
$(\mathbb{Z}_2)^{n-1}\!\rtimes S_n$.

The groups of signed and even-signed permutations of~$O_n$ are isomorphic
to the Coxeter groups $\mathcal{B}_n$ (when $n\ge2$) and $\mathcal{D}_n$
(when $n\ge4$), respectively.  This makes possible the use of known results
on subgroups, cosets, generators, numbers of elements with given
order,~etc.  The following is a brief explanation of these isomorphisms
(cf.~\cite{Brenti94}).  A~Coxeter group with generators
$\sigma_1,\dots,\sigma_n$ is defined by the relations $\sigma_i^2=1$ for
each~$i$, and zero or more `braid relations' of the form
${(\sigma_i\sigma_j)^{c_{ij}}=1}$.  These may be encoded as a Coxeter graph
on the vertex set $\sigma_1,\dots,\sigma_n$, with each
pair~$\sigma_i,\sigma_j$ joined by $c_{ij}-2$~edges.  The non-isomorphic
Coxeter groups $\mathcal{A}_n$ ($n\ge1$), $\mathcal{B}_n$ ($n\ge2$), and
$\mathcal{D}_n$ ($n\ge4$) are the groups defined respectively by
(1)~$c_{i,i+1}=3$ for each~$i\le n-1$; (2)~$c_{i,i+1}=3$ for each~$i\le
n-2$, and $c_{n-1,n}=4$; and (3)~$c_{i,i+1}=3$ for each~$i\le n-2$, and
$c_{n-2,n}=3$.  The group~$\mathcal{A}_n$ turns~out to be isomorphic to the
symmetric group~$S_{n+1}$.  It~follows from the associated graphs that
$\mathcal{D}_3\cong\mathcal{A}_3\cong S_4$; the notation~$\mathcal{D}_3$ is
not standard.

That the group of signed permutations of the $n$-set~$O_n$ is isomorphic
to~$\mathcal{B}_n$ (when~$n\ge2$) follows by defining
$\sigma_i=[i_+(i+1)_+]$ for $i=1,\dots,n-1$, and $\sigma_n=[n_-]$.
Similarly, that the group of even-signed permutations is isomorphic
to~$\mathcal{D}_n$ (when~$n\ge4$) follows by defining
$\sigma_1,\dots,\sigma_{n-1}$ in the same way, and
$\sigma_n=[({n-1})_-n_-]$.  Each of the cases $n=1,2,3$ is at~least partly
degenerate.  The groups $\mathcal{B}_n,\mathcal{D}_n$ are respectively
isomorphic, when $n=1$, to $\mathbb{Z}_2$ and the trivial group;
when~$n=2$, to the order-$8$ dihedral group and the Klein $4$-group
$\mathbb{Z}_2\times\mathbb{Z}_2$; and when~$n=3$, to an order-$48$ group
and~$\mathcal{A}_3\cong S_4$.  It~is easy to see that when $n$~is odd,
$\mathcal{B}_n$~is isomorphic to~$\mathbb{Z}_2\times\mathcal{D}_n$, where
the $\mathbb{Z}_2$ factor is generated by $[1_-]\cdots[n_-]$.  So
$\mathcal{B}_3$~is isomorphic to $\mathbb{Z}_2\times S_4$, though this fact
will not be used in what follows.

\section{Point Configurations}
\label{sec:configs}

To study the effects of M\"obius transformations on the reduced forms
(\ref{eq:fuchsian0}),(\ref{eq:fuchsian1}), and their interplay with index
transformations, it~is helpful to focus first on their effects on the
singular points.  The action of the M\"obius group $PGL(2,\mathbb{C})$ on
the complex projective line~$\mathbb{P}^1(\mathbb{C})$ is {\em sharply
$3$-transitive\/}: for any two triples of distinct points
$(x_1,x_2,x_3),\allowbreak(x_1',x_2',x_3')$, there is a unique $P\in
PGL(2,\mathbb{C})$ satisfying $x_i'=Px_i$, $i=1,2,3$.  Due~to this, it~is
traditional to assume that the singular points of Fuchsian differential
equations with $n\ge3$ include~$0,1,\infty$, thereby quotienting~out
(partially) the action of~$PGL(2,\mathbb{C})$.

If $X=\{x_1,\dots,x_n\}$ is an unordered set of $n$~distinct points, three
of which are~$0,1,\infty$, there will be exactly $n(n-1)(n-2)$
transformations $P\in PGL(2,\mathbb{C})$ for which $PX$~also has this
property.  For computational purposes, this is best restated in a form
referring to {\em ordered\/} sets of points.

\begin{definition}
\mbox{}
\begin{enumerate}
\item 
\begin{enumerate}
\item $X[n]\defeq[\mathbb{P}^1(\mathbb{C})]^n\setminus\Delta$ is the space
of $n$~distinct (labeled) points in~$\mathbb{P}^1(\mathbb{C})$, where
$\Delta$~is the set of configurations with two points coinciding.  The
group $PGL(2,\mathbb{C})$ acts on~$X[n]$ diagonally, i.e., on each point
individually, and $S_n$~acts on~it by permuting points.  Since these
actions commute, the direct product $PGL(2,\mathbb{C})\times S_n$ also acts
on~$X[n]$.
\item When $n\ge3$ and $1\le i_1,i_2,i_3\le n$ are distinct,
$X_{i_1,i_2,i_3}[n]\subset X[n]$ is the space of $n$-tuples satisfying
$x_{i_1}=0$, $x_{i_2}=1$, $x_{i_3}=\infty$.
\end{enumerate}
\item $X(n)\defeq X[n]/PGL(2,\mathbb{C})$ is the space of $n$~distinct
(labeled) points in $\mathbb{P}^1(\mathbb{C})$, up~to M\"obius
transformations.  When~$n\ge3$, $X(n)$~may be realized as any of
the~$X_{i_1,i_2,i_3}[n]$, in particular as $X_{1,2,3}[n]$.
\end{enumerate}
\end{definition}

\begin{remarkaftertheorem}
The relation between $n$-tuples $x_{i_1,i_2,i_3}=(x_1,\dots,x_n)\in
X_{i_1,i_2,i_3}[n]$ and $x_{i'_1,i'_2,i'_3}=(x'_1,\dots,x'_n)\in
X_{i_1',i_2',i_3'}[n]$, both representing the same point in~$X(n)$, is
straightforward: $x'_i=Px_i$ for each~$i$, where $P\in PGL(2,\mathbb{C})$
is the unique M\"obius transformation that takes
$x_{i'_1},x_{i'_2},x_{i'_3}$ to~$0,1,\infty$.
\end{remarkaftertheorem}

\begin{definition}
For every $\mathfrak{x}\in X_{1,2,3}[n]\subset X[n]$, the associated
$n$-point hyper-Kummer group $\mathfrak{K}_n(\mathfrak{x})$ is the subset
of the direct product group $PGL(2,\mathbb{C})\times S_n$, acting
on~$X[n]$, which contains all elements~$g$ for which $g\cdot
\mathfrak{x}\in X_{1,2,3}[n]$.
\end{definition}

\begin{proposition}
\label{prop:1}
For every $\mathfrak{x}\in X_{1,2,3}[n]$, $\mathfrak{K}_n(\mathfrak{x})$~is
a group isomorphic to~$S_n$, and the collection of\/ $n!$~elements
$(P,s)\in PGL(2,\mathbb{C})\times S_n$ which $\mathfrak{K}_n(\mathfrak{x})$
contains can be generated as follows.  Let
$\mathfrak{x}=(x_1,\dots,x_n)=(0,1,\infty,a_1,\dots,a_{n-3})$.  For each
choice of distinct\/ ${1\le i_1,i_2,i_3\le n}$, there is a unique\/ $P\in
PGL(2,\mathbb{C})$ that takes\/ $x_{i_1},x_{i_2},x_{i_3}$
to\/~$0,1,\infty$.  There are\/ $n(n-1)(n-2)$ such choices.  For each
such~$P$, choose $s\in S_n$, a~permutation of\/~$1,\dots,n$, such that
$(Px_{s(1)},\dots,Px_{s(n)})\in X_{1,2,3}[n]$.  Since necessarily
$s(k)=i_k$, $k=1,2,3$, there are $(n-3)!$ choices for~$s$.
\end{proposition}

Informally, any $\sigma\in\mathfrak{K}_n(\mathfrak{x})\cong S_n$, where
$\mathfrak{x}=(0,1,\infty,a_1,\dots,a_{n-3})$, acts as the composition of
two transformations: (i)~the M\"obius transformation which, acting
diagonally, takes three of these $n$~points to $0,1,\infty$ respectively,
and (ii)~a permutation of the images of the remaining $n-3$~points under
this transformation.

\begin{example}
\label{ex:hyperg}
$n=3$, the hypergeometric case.  This case is quite degenerate:
$X_{1,2,3}[3]=\{(0,1,\infty)\}$ is a single point, as is~$X(3)$.  The
$3$-point hyper-Kummer group,
$\mathfrak{K}_3((0,1,\infty))<PGL(2,\mathbb{C})\times S_3$, is the
traditional Kummer group~\cite{Prosser94}, also called the cross-ratio
group.  It~is based on the six M\"obius transformations $x\mapsto f(x)$,
with~$f(x)$ equal to
\begin{equation}
\label{eq:kummergroup}
x,\qquad \frac x{x-1},\qquad 1-x,\qquad \frac{x-1}x,\qquad \frac1x,\qquad \frac1{1-x}.
\end{equation}
These permute $0,1,\infty$, so each must be followed by an appropriate
(unique) $s\in S_3$, to yield a transformation in $PGL(2,\mathbb{C})\times
S_3$ that takes $(0,1,\infty)$ to itself.
\end{example}

\begin{example}
\label{ex:heun}
$n=4$, the Heun case.  $X_{1,2,3}[4]$ is the space of configurations of the
form $(0,1,\infty,a)$, so it is effectively the triply punctured sphere
$\mathbb{P}^1(\mathbb{C})\setminus\{0,1,\infty\}\ni a$, as is~$X(4)$.  The
$4$-point hyper-Kummer group, isomorphic to~$S_4$, was worked~out by Heun
and may be called the Heun group.  $\mathfrak{K}_4((0,1,\infty,a))$~is
based on the $24$~M\"obius transformations $x\mapsto f(x)$, with $f(x)$
equal to
\begin{equation}
\label{eq:heungroup}
\begin{array}{cccccc}
x, & \dfrac{x}{x-1}, & \dfrac{x}{a}, & \dfrac{x}{x-a}, & \dfrac{(1-a)x}{x-a}, & \dfrac{(a-1)x}{a(x-1)},\\[12pt]
1-x, & \dfrac{x-1}{x}, & \dfrac{x-1}{x-a}, & \dfrac{x-1}{a-1}, &
\dfrac{a(x-1)}{x-a}, & \dfrac{a(x-1)}{(a-1)x},\\[12pt]
\dfrac{1}x, & \dfrac{1}{1-x}, & \dfrac{a}{x}, & \dfrac{a}{a-x}, &
\dfrac{a-1}{x-1}, & \dfrac{1-a}{x-a},
\\[12pt]
\dfrac{x-a}x, & \dfrac{a-x}a, & \dfrac{x-a}{x-1}, & \dfrac{a-x}{a-1}, & \dfrac{x-a}{a(x-1)}, & \dfrac{a-x}{(a-1)x}.
\end{array}
\end{equation}
The transformations in the four rows of~(\ref{eq:heungroup}) take
$0,1,\infty,a$ respectively to~$0$, and in a sense, each permutes
$0,1,\infty,a$.  What this means is that for any $(0,1,a,\infty)\in
X_{1,2,3}[4]$, i.e., any
$a\in\mathbb{P}^1(\mathbb{C})\setminus\{0,1,\infty\}$, each takes the {\em
unordered\/} set $\{0,1,\infty,a\}$ to a set of the form
$\{0,1,\infty,a'\}$.  An example is $x\mapsto\frac{x}{x-a}$, which takes
$0,1,\infty,a$ to $0,a'=\frac1{1-a},1,\infty$, respectively.  This element
of $PGL(2,\mathbb{C})$ can be written in cycle notation as~$(0)(1a\infty)$,
with the understanding that ``$1$~is mapped to~$a$'' really means ``$1$~is
mapped to~$a'$\,''.  It~must be followed by an appropriate (unique)
permutation $s\in S_4$, if it is to take $(0,1,\infty,a)\in X_{1,2,3}[4]$
to $(0,1,\infty,a')\in X_{1,2,3}[4]$.

The $a$-dependent transformation groups
$\mathfrak{K}_4((0,1,\infty,a))\subset PGL(2,\mathbb{C})\times S_4$, each
isomorphic to~$S_4$, can usefully be viewed as arising from the action of a
{\em single\/} $4$-point hyper-Kummer group isomorphic to~$S_4$, which may
be denoted~$\mathfrak{K}_4$, on~$X_{1,2,3}[4]$.  This~is because each
formal permutation of~$0,1,\infty,a$ gives a map from
$a\in\mathbb{P}^1(\mathbb{C})\setminus\{0,1,\infty\}$ to
$a'\in\mathbb{P}^1(\mathbb{C})\setminus\{0,1,\infty\}$, i.e., from
$X_{1,2,3}[4]$ to itself.  A~similar interpretation is possible for
any~$n\ge3$, but the $n=4$ case is slightly degenerate, though not so much
so as the $n=3$ case.  By examination, the representation of~$S_4$ by
transformations of~$X_{1,2,3}[4]$ is not faithful.  The orbit of any
$a\in\mathbb{P}^1(\mathbb{C})\setminus\{0,1,\infty\}$ consists of no~more
than $6$~points $a'\in\mathbb{P}^1(\mathbb{C})\setminus\{0,1,\infty\}$,
rather than~$24$.  These can be obtained directly from~$a$ by $a\mapsto
a'=f(a)$, with $f$~ranging over the M\"obius transformations
of~(\ref{eq:kummergroup}).  The reduction of $24$~values to~$6$ arises from
the invariance of any~$a$ under an order-$4$ normal subgroup
of~$\mathfrak{K}_4\cong S_4$ isomorphic to the Klein $4$-group
$\mathbb{Z}_2\times\mathbb{Z}_2$, which comprises the maps $x\mapsto x,\
\frac{a}x,\ \frac{x-a}{x-1},\ \frac{a(x-1)}{x-a}$, or in cycle notation,
$(0)(1)(\infty)$, $(0\infty)(1a)$, $(0a)(1\infty)$, $(01)(a\infty)$.  For
each, $a'=a$.

Generic $a$-orbits contain $6$~distinct values~$a'$, but two nongeneric
ones contain fewer: the harmonic (lemniscatic) orbit $\{-1,\frac12,2\}$ and
the equianharmonic orbit~$\{\frac12\pm{\rm i}\frac{\sqrt3}2\}$.  The names
are taken from elliptic function theory~\cite[chap.~18]{Abramowitz64}.
These small orbits contain values of~$a$ that are invariant under larger
subgroups of~$S_4$ than the Klein $4$-group.  Each of $a=-1,\frac12,2$ is
fixed by a dihedral group of order~$8$, the three groups being conjugate
in~$S_4$.  Similarly, each of $a=\frac12\pm{\rm i}\frac{\sqrt3}2$ is fixed
by the alternating subgroup~$A_4$.
\end{example}

\begin{example}
$n=5$, a nonclassical case.  $X_{1,2,3}[5]$ is the space of $5$-point
configurations in~$\mathbb{P}^1(\mathbb{C})$ of the form
$(0,1,\infty,a,b)$, so it is effectively the complement of five complex
lines in~$\mathbb{C}^2$, i.e., $\{(a,b)\in\mathbb{C}^2\mid a\neq0,1,\
b\neq0,1,\ a\neq b\}$, as is~$X(5)$.  However, the action on~$X_{1,2,3}[5]$
of the $5$-point hyper-Kummer group~$\mathfrak{K}_5$, isomorphic to~$S_5$,
is difficult to visualize.  The $5!=120$ transformations in this group are
too numerous to list, but the following example should clarify their
effects on~$(a,b)$.

By Proposition~\ref{prop:1}, the action of each $\sigma\in
\mathfrak{K}_5\cong S_5$ on~$(0,1,\infty,a,b)$ is based on a M\"obius
transformation that maps three of $0,1,\infty,a,b$ to~$0,1,\infty$.
Suppose $a,b,0$ are mapped to $0,1,\infty$, respectively.  Then this
transformation must be $x\mapsto \frac{b(x-a)}{(b-a)x}$.  The images of the
remaining two points, namely $1,\infty$, will be the transformed
coordinates $a',b'$; in~either order, since the action of~$\sigma$ involves
a permutation of the remaining points, in~addition to the M\"obius
tranformation.  It~follows that the action of~$\sigma$ on~$X_{1,2,3}[5]$
will be either the map $(a,b)\mapsto(a',b')$ with
\begin{equation}
a'=\frac{b(1-a)}{b-a},\qquad b'=\frac{b}{b-a},
\end{equation}
or the same with $a',b'$ interchanged.  The corresponding
$\sigma\in\mathfrak{K}_5\cong S_5$ would be written in cycle notation as
$(a0\infty b1)$ and~$(a0\infty)(1b)$.  The above caveat applies: for
instance, in the first, ``$\infty$~is mapped to~$b$'' really means
``$\infty$~is mapped to~$b'$\,''.

The $n=5$ case is not degenerate, unlike $n=3,4$. If $(0,1,\infty,a,b)\in
X_{1,2,3}[5]$ is `generic', i.e., if
$0,1,\infty,a,b\in\mathbb{P}^1(\mathbb{C})$ are in general position, then
its orbit under $\mathfrak{K}_5\cong S_5$ will consist of $120$~distinct
points of the form $(0,1,\infty,a',b')$.  Such an orbit is effectively a
set of $120$~distinct points $(a',b')$ in the complement of $5$~lines
in~$\mathbb{C}^2$, mentioned above. But there are nongeneric orbits.  Any
orbit that includes a point $(a,b)$ with $a+b=1$ will contain no~more than
$60$~points, since such a point is left invariant by a $\mathbb{Z}_2$-group
comprising the maps $(a,b)\mapsto(a,b)$, $(a,b)\mapsto(1-b,1-a)$, or in
cycle notation, $(0)(1)(\infty)(a)(b)$, $(01)(ab)(\infty)$.  Each point
$(a',b')$ on such an orbit will be invariant under a $\mathbb{Z}_2$-group
which is conjugate in~$S_5$ to this one.

Some nongeneric orbits of this type are even smaller.  It~can be shown that
the orbits of $(a,b)=(\frac{1+i}2,\frac{1-i}2)$, $(\frac{1+{\rm
i}{\sqrt3}}2,\frac{1-{\rm i}{\sqrt3}}2)$,
$(\frac{3+\sqrt5}2,\frac{-1-\sqrt5}2)$ consist of $30,20,12$ points
$(a',b')$, respectively.
\end{example}

\section{Key Results}
\label{sec:key}

The joint effects of M\"obius and index transformations on Fuchsian
differential equations with $n$~singular points can now be determined.
To~maintain contact with the classical theory of special functions, three
of the singular points will be fixed at~$x=0,1,\infty$.  That~is, the
finite singular points $d_1,\dots,d_{n-1}$ will be taken to be
$0,1,a_1,\dots,a_{n-3}$.  The results of~\S\,\ref{sec:configs} on how a
group $\mathfrak{K}_n\cong S_n$~acts on ordered point configurations of the
form $(0,1,\infty,a_1,\dots,a_{n-3})$ will then apply.

\begin{definition}
The collection of joint M\"obius--index transformations that leave
invariant, up~to parameter changes, a version of the symmetrically reduced
form~(\ref{eq:fuchsian0}) that includes $x=0,1$ among the finite singular
points, denoted
\begin{equation}
\label{eq:sym}
\mathcal{E}_{\rm sym}(0,\delta^{(0)};\,1,\delta^{(1)};\,\infty,\delta^{(\infty)};\,a_1,\delta_1;\dots;a_{n-3},\delta_{n-3};\,Q),
\end{equation}
and a similar restricted version of the asymmetrically reduced
form~(\ref{eq:fuchsian1}), denoted
\begin{equation}
\label{eq:asym}
\mathcal{E}_{\rm asym}(0,\rho^{(0)};\,1,\rho^{(1)};\,\infty,\rho^{(\infty)}\!,{\hat\rho}^{(\infty)};\,a_1,\rho_1;\dots;a_{n-3},\rho_{n-3};\,Q),
\end{equation}
will be called the automorphism group of symmetrically, resp.\
asymmetrically, reduced $n$-point Fuchsian differential equations.  By
convention, permutations of $(a_1,\delta_1),\dots,(a_{n-3},\delta_{n-3})$,
resp.\ $(a_1,\rho_1),\dots,(a_{n-3},\rho_{n-3})$, are included, though
interchange of the exponents $\rho^{(\infty)}\!,{\hat\rho}^{(\infty)}$ at
infinity is~not.
\end{definition}

Each of these groups acts on the associated parameter space, which is
coordinatized by the $n-3$ singular point locations $a_1,\dots,a_{n-3}$,
$n$~(independent) exponent parameters, and $n-3$ accessory parameters.

\begin{theorem}
\label{thm:symmetric}
The automorphism group $G_{\rm sym}(n)$ of symmetrically reduced $n$-point
Fuchsian differential equations is isomorphic to~$S_n$.
\end{theorem}
\begin{proof}
If $\mathfrak{x}=(0,1,\infty,a_1,\dots,a_{n-3})$ is the set of singular
points, any $\sigma\in\mathfrak{K}_n(\mathfrak{x})\cong S_n$ consists of a
M\"obius transformation~$P_\sigma$ that acts on~$\mathfrak{x}$ diagonally,
and a subsequent permutation.  But under M\"obius transformations of the
independent variable, characteristic exponents of Fuchsian equations
accompany the points at~which they are evaluated.  This follows from the
general theory of differential equations in the complex domain, formalized
using the Riemann P-symbol \cite[\S\,15.3]{Erdelyi53}.  So under~$\sigma$,
(\ref{eq:sym})~will typically be transformed to another~$\mathcal{E}_{\rm
sym}$, with transformed parameters.  One case requires special treatment:
when $d=P_\sigma\infty\neq\infty$.  In~this case, the sum of the two
exponents at the finite singular point $x'=d$ and at~$x'=\infty$ will be
$-1$~and~$1$, respectively, rather than the other way around, as they
should~be.  So $\sigma$~must be followed by a uniquely determined index
transformation, based on $S=S(x')=x'-d$, which will shift the exponents,
and yield a proper~$\mathcal{E}_{\rm sym}$.
\end{proof}

\begin{theorem}
\label{thm:asymmetric}
The automorphism group $G_{\rm asym}(n)$ of asymmetrically reduced
$n$-point Fuchsian differential equations is isomorphic to the Coxeter
group\/~$\mathcal{D}_n$, the group of even-signed permutations of an
$n$-set, i.e., is isomorphic to a semi-direct product of the form\/
$(\mathbb{Z}_2)^{n-1}\!\rtimes S_n$.
\end{theorem}
\begin{remarkaftertheorem}
The subgroup $S_n$ here is $\mathfrak{K}_n$, i.e., the group of M\"obius
automorphisms, including subsequent permutations of the singular points,
and the normal subgroup $(\mathbb{Z}_2)^{n-1}$ is the group of
index-transformation automorphisms defined in~\S\,\ref{sec:basic}.
In~cycle notation they comprise, respectively, all unsigned
(or~equivalently positively signed) permutations
of~$0,1,\infty,a_1,\dots,a_{n-3}$, and all even-signed permutations of the
form $[0_\pm][1_\pm][\infty_\pm][(a_1)_\pm]\dotsb[(a_{n-3})_\pm]$.
\end{remarkaftertheorem}
\begin{proof}
Each $\sigma\in\mathfrak{K}_n\cong S_n$ yields an automorphism
of~(\ref{eq:asym}), i.e., of~(\ref{eq:fuchsian1}), which transforms the
variable~$x$ to~$P_\sigma x$.  If $\sigma(\infty)\neq\infty$, the
transformed equation's exponents at $x'=d=\sigma(\infty)$ will be
$\rho^{(\infty)},{\hat\rho}^{(\infty)}$, and must be shifted so that one of
them is zero, to preserve the form of~(\ref{eq:asym}).  By~convention,
$\rho^{(\infty)}$~will be chosen.  The index transformation based on
$S=S(x')=(x'-d)^{-\rho^{(\infty)}}$ \!\!will perform this shift.

As a set, the collection of joint M\"obius--index transformations $G_{\rm
asym}(n)$ that leave~(\ref{eq:fuchsian1}) invariant will be the Cartesian
product $(\mathbb{Z}_2)^{n-1}\!\times S_n$, since the group of
index-transformation automorphisms is isomorphic to $(\mathbb{Z}_2)^{n-1}$.
(Any index transformation involves a binary choice at each of the $n-1$
finite singular points.)  If $(h,\sigma)\in G_{\rm asym}(n)$, the index
transformation~$h$ will be viewed as acting {\em after\/}~$\sigma\in S_n$.

The structure of~$G_{\rm asym}(n)$ {\em as a group\/} remains to be
determined.  In~resolving this issue, a seemingly minor point will be
important.  The two exponents $\rho^{(\infty)}\!,{\hat\rho}^{(\infty)}$
do~not appear explicitly in the differential equation~(\ref{eq:fuchsian1}),
though they~do in the symbolic form~(\ref{eq:asym}).  They are uniquely
determined by~(\ref{eq:fuchsian1}) only up~to interchange.  When
determining the group structure of~$G_{\rm asym}(n)$, a~convention for
ordering this pair, after applying any $\sigma\in S_n$ and
$h\in(\mathbb{Z}_2)^{n-1}$, will need to be followed.

Any automorphism $(h,\sigma)\in G_{\rm asym}(n)$ may be written as a signed
permutation of the singular points~$(0,1,\infty,a_1,\dots,a_{n-3})$, i.e.,
as some $g\in\mathcal{B}_n$, by annotating the cycle representation
of~$\sigma$ as~follows.  If $d$~is any of the $n$~points, then provided
$\sigma(d)\neq\infty$, the annotation $\dots d_+\sigma(d)\dots$, resp.\
$\dots d_-\sigma(d)\dots$, will indicate that $d\mapsto\sigma(d)$ is,
resp.\ is~not, followed by a negation of the nonzero exponent
at~$x'=\sigma(d)$.  When $\sigma(d)=\infty$, the annotation of $\dots
d\infty\dots$ will have a special meaning.  The annotation $\dots
d_-\infty\dots$, as opposed to $\dots d_+\infty\dots$, will signify that
after the index transformation~$h$ is performed, the exponent parameters
$\rho^{(\infty)}\!,{\hat\rho}^{(\infty)}$~are interchanged.  Such a
transposition will be included, if necessary, to ensure {\em
even-signedness\/}: an~even number of negative signs in the annotated cycle
representation of~$g$.  With this convention, $g\in
\mathcal{D}_n\!\vartriangleleft \mathcal{B}_n$.

What remains to be shown is that the bijection between $G_{\rm asym}(n)$
and~$\mathcal{D}_n$ is a group isomorphism.  The fact that
$\mathcal{D}_n\cong (\mathbb{Z}_2)^{n-1}\!\rtimes_\Gamma S_n$ can be used,
where $\Gamma:S_n\to{\rm Aut}((\mathbb{Z}_2)^{n-1})$ denotes the action
of~$S_n$ by conjugation.  If $\sigma\in S_n$ fixes~$\infty$, i.e., permutes
only the $n-1$ finite singular points, then conjugating
$h\in(\mathbb{Z}_2)^{n-1}$ by~$\sigma$ will permute the points at~which
$h$~negates an exponent.  So in this case, the conjugation actions on the
copies of~$(\mathbb{Z}_2)^{n-1}$ in $\mathcal{D}_n$ and~$G_{\rm asym}(n)$
are isomorphic.

The harder case is when $\sigma(\infty)\neq\infty$.  There are several
subcases, but the following example makes clear how each can be handled.
Let $\sigma$ be the involution~$(0\infty)$, or in annotated notation
$[0_+\infty_+]$, and let $h=[1_-][\infty_-]$.  The equation
\begin{equation}
\label{eq:signedpermproduct}
[0_+\infty_+]^{-1} \cdot [1_-][\infty_-] \cdot
[0_+\infty_+] = [0_-][1_-][\infty_+]
\end{equation}
holds in~$\mathcal{D}_n$, but its validity in~$G_{\rm asym}(n)$ is not
obvious.  The automorphism labeled by $[1_-][\infty_-]$ negates the
exponent parameter at {\em one\/} of the $n-1$~finite singular points, but
the conjugated automorphism labeled by $[0_-][1_-][\infty_+]$ does so
at~{\em two\/}.  The identity of the left side with the right can be
verified by sequentially working~out the action of the product on the left.
Take $G_{\rm asym}(n)$ to act on~(\ref{eq:asym}), i.e.,
on~$\mathcal{E}_{\rm asym}$, `on~the left', so that the order of action of
group elements, not previously specified, will be from right to left.  Then
some straightforward computations reveal that 
{
\setlength\multlinegap{0.65in}
\begin{enumerate}
\item $[0_+\infty_+]\mathcal{E}_{\rm asym}$, obtained by the M\"obius
transformation $x\mapsto x'=1/x$ from~$\mathcal{E}_{\rm asym}$, is of
the form
\begin{multline*}
\mathcal{E}_{\rm asym}(0,{\hat\rho}^{(\infty)}\!-\rho^{(\infty)};\,1,\rho^{(1)};\,\infty,\rho^{(\infty)}\!,\rho^{(0)}\!+\rho^{(\infty)};\\
\setlength\multlinegap{0.3in}
\hfill\tfrac1{a_1},\rho_1;\dots;\tfrac1{a_{n-3}},\rho_{n-3};\,Q_1),
\end{multline*}
with independent variable $x$;
\item $[1_-][\infty_-] \cdot[0_+\infty_+]\mathcal{E}_{\rm asym}$, obtained
by a further index transformation at $x'=1$, is of the form
\begin{multline*}
\mathcal{E}_{\rm asym}(0,{\hat\rho}^{(\infty)}\!-\rho^{(\infty)};\,1,-\rho^{(1)};\,\infty,\rho^{(0)}\!+\rho^{(1)}\!+\rho^{(\infty)}\!,\rho^{(1)}\!+\rho^{(\infty)};\\
\setlength\multlinegap{0.3in}
\tfrac1{a_1},\rho_1;\dots;\tfrac1{a_{n-3}},\rho_{n-3};\,Q_2),
\end{multline*}
with independent variable $x'$; and
\item $[0_+\infty_+]^{-1} \cdot [1_-][\infty_-]
\cdot[0_+\infty_+]\mathcal{E}_{\rm asym}$, obtained by a further
M\"obius transformation $x'\mapsto x=1/x'$, is of the form
\begin{multline*}
\mathcal{E}_{\rm asym}(0,-\rho^{(0)};\,1,-\rho^{(1)};\,\infty,\rho^{(0)}\!+\rho^{(1)}\!+\rho^{(\infty)}\!,\rho^{(0)}\!+\rho^{(1)}\!+{\hat\rho}^{(\infty)};\\
\setlength\multlinegap{0.3in}
a_1,\rho_1;\dots;a_{n-3},\rho_{n-3};\,Q_3),
\end{multline*}
with independent variable~$x$.
\end{enumerate}
}
\noindent
The last of these, with negated exponents at $x=0,1$, clearly has the same
exponent parameters as $[0_-][1_-][\infty_+]\mathcal{E}_{\rm asym}$, i.e.,
as (\ref{eq:asym})~acted~on by the right side
of~(\ref{eq:signedpermproduct}).  The vectors of accessory parameters can
also be shown to correspond, by a computation omitted here.  This careful
sequential approach can be applied in~general, to show that conjugation of
$(\mathbb{Z}_2)^{n-1}$ \!by~$\sigma\in S_n$ always has the same meaning
in~$G_{\rm asym}(n)$ as in~$\mathcal{D}_n$.
\end{proof}

Asymmetrically reduced Fuchsian differential equations predominate in
applied mathematics, though symmetrically reduced ones are used in
conformal mapping.  It~is usual to express solutions of~(\ref{eq:asym})
in~terms of a {\em canonical local solution\/}: the Frobenius solution
at~$x=0$ that belongs to the exponent zero, rather than the other
exponent~$\rho^{(0)}$.  It~will be denoted
\begin{equation}
F(0,\rho^{(0)};\,1,\rho^{(1)};\,\infty,\rho^{(\infty)}\!,{\hat\rho}^{(\infty)};\,a_1,\rho_1;\dots;a_{n-3},\rho_{n-3};\,Q;\,x).
\end{equation}
If $\rho^{(0)}$~is not a positive integer, $F$~is guaranteed to be analytic
in a neighborhood of~$x=0$, and may be chosen to equal unity at~$x=0$.

\begin{proposition}
\label{prop:generate}
A family of\/ $2^{n-1}n!$ formally distinct local solutions
of\/~{\rm(\ref{eq:asym})}, indexed by the group $G_{\rm asym}(n)\cong
\mathcal{D}_n$, can be generated as~follows.  Any $g\in G_{\rm asym}(n)$
will include a M\"obius transformation $x\mapsto x'=Px$, together with a
subsequent permutation of singular points and a subsequent index
transformation.  Apply $g$ to~{\rm(\ref{eq:asym})}, and compute the
zero-exponent Frobenius solution of the transformed equation at its
singular point~$x'=0$.  This will be of the form $A(x)$ times
\begin{displaymath}
F(0,{\rho}^{(0)\prime};\,1,{\rho}^{(1)\prime};\,\infty,{\rho}^{(\infty)\prime}\!,{\hat\rho}^{(\infty)\prime};\,a'_1,\rho'_1;\dots;a'_{n-3},\rho'_{n-3};\,Q';\,Px),
\end{displaymath}
where
\begin{displaymath}
A(x)={\rm const}\times x^{\nu^{(0)}}(x-1)^{\nu^{(1)}}\prod_{i=1}^{n-3}
(x-a_i)^{\nu_i},
\end{displaymath}
and where
$\rho^{(0)\prime}\!,\rho^{(1)\prime}\!,\rho^{(\infty)\prime}\!,{\hat\rho}^{(\infty)\prime}\!,\rho'_1,\dots,\rho'_{n-3}$
and also $\nu^{(0)},\nu^{(1)},\nu_{1},\dots,\nu_{{n-3}}$ are linear
functions of\/
$\rho^{(0)},\rho^{(1)},\rho^{(\infty)},{\hat\rho}^{(\infty)},\rho_1,\dots,\rho_{n-3}$,
with coefficients in\/~$\{0,\pm1\}$.
\end{proposition}
\begin{proof}
Any $g\in G_{\rm asym}(n)$ is a pair $(h,\sigma)\in
(\mathbb{Z}_2)^{n-1}\!\rtimes S_n$, with $\sigma$~acting first.  The
automorphism~$\sigma$ changes the independent variable from $x$
to~$x'=P_\sigma x$, where $P_\sigma$~is the associated M\"obius
transformation.  The index transformation~$h$ acts by replacing the
differential operator~$T_1$ in the equation by $\widehat T_1=ST_1S^{-1}$,
where $S=S(x')$ is a product over finite transformed singular points of the
form $\prod_{i=1}^{n-1}(x'-d_i')^{s_i'}$, with $s_i'$~the amount by which
the exponents at~$x'=d_i'$ are shifted.  As~was noted in the proof of
Theorem~\ref{thm:asymmetric}, if $\sigma(\infty)\neq\infty$ then the
automorphism~$\sigma$ will contribute an additional factor of the form
$(x'-d_i')^{s_i'}$ to the conjugating function~$S$.  Conjugating $T_1$ by
the resulting~$S$ is equivalent to replacing any solution $u=u(x')$ of the
transformed differential equation, such as $u(x')=F(x')=F(P_\sigma x)$, by
$S(x')u(x')$.  Up~to a constant factor, $S$~can be written as a product
over the {\em untransformed\/} finite singular points,
$\prod_{i=1}^{n-1}(x-d_i)^{s_i}$, since the transformed singular points are
the original singular points, acted on by~$P_\sigma$ and possibly permuted.
\end{proof}

\begin{theorem}
\label{thm:local}
Suppose $g\in G_{{\rm asym},0}(n)\cong
\mathcal{D}_{n-1}\cong(\mathbb{Z}_2)^{n-2}\!\rtimes S_{n-1}$, the subgroup
of~$G_{\rm asym}(n)$ that comprises all automorphisms that fix $x=0$ and
perform no~index transformation there. That~is, in annotated cycle notation
$g$~should contain the positive\/ $1$-cycle\/~$[0_+]$.  Then the local
solution of\/~{\rm(\ref{eq:asym})} obtained from~$g$ by the technique of
Proposition\/~{\rm\ref{prop:generate}}, as a function of~$x$ near~$x=0$,
will equal $F$~itself, pointwise.
\end{theorem}
\begin{proof}
If $g=(h,\sigma)$ with the M\"obius transformation~$P_\sigma$ taking $x=0$
to~$x'=0$, then $F(Px)$~will be analytic at~$x=0$, and moreover the
prefactor~$A(x)$ will not include an $x^{\nu^{(0)}}$ \!\!factor.  So it
will belong to the zero exponent at~$x=0$, and must be a multiple of $F$
itself.
\end{proof}

\begin{corollary}
\label{cor:transformationgroup}
Of the $2^{n-1}n!$ formally distinct local solutions
of~{\rm(\ref{eq:asym})} generated by the technique of
Proposition\/~{\rm\ref{prop:generate}}, $2^{n-2}(n-1)!$ solutions\/
{\rm(}including~$F${\rm)} are equivalent expressions for~$F$; and since
this collection is bijective with\/~$\mathcal{D}_{n-1}$, the family of
transformations of~$F$ into any of these alternative forms, which will be
called the\/ {\rm transformation group of}~$F$, has a group structure
isomorphic to\/~$\mathcal{D}_{n-1}$.
\end{corollary}

\begin{theorem}
The\/ $2^{n-1}n!$ local solutions of\/~{\rm(\ref{eq:asym})} indexed by
$G_{{\rm asym}}(n)\cong \mathcal{D}_{n}$ split into\/ $2n$~sets of\/
$2^{n-2}(n-1)!$~equivalent expressions involving~$F$, each set defining one
of the two Frobenius solutions at one of the $n$~singular points.  These\/
$2n$ sets are bijective with the cosets of $G_{{\rm asym},0}(n)\cong
\mathcal{D}_{n-1}$ in $G_{{\rm asym}}(n)\cong\mathcal{D}_n$.
\end{theorem}
\begin{proof}
This is similar to the proof of Theorem~\ref{thm:local}.  For any of these
$2n$~sets, the identification of every solution in the set with a single
Frobenius solution follows from their identical local behavior at the
defining singular point.
\end{proof}

A normalization of the $2^{n-1}n!$ formally distinct local solutions
of~(\ref{eq:asym}) has not yet been chosen.  A~reasonable one would take
each solution defined near a finite singular point~$x=d$ with
exponent~$\rho$ to be asymptotic to $(x-d)^\rho$ or $(d-x)^\rho$ as~$x\to
d$.  Similarly, those with exponents
$\rho^{(\infty)}\!,{\hat\rho}^{(\infty)}$ at~$x=\infty$ would be taken to
be asymptotic to $x^{-\rho^{(\infty)}}\!\!,x^{-{\hat\rho}^{(\infty)}}$
\!\!as~$x\to\infty$.  These conventions assume the exponent difference is
not an integer, so that both solutions display power-law behavior.

The following discussion should clarify the action of $G_{\rm
asym}(n)\cong\mathcal{D}_n$ on the asymmetrically reduced equation, and in
particular, the action of~$\mathcal{D}_{n-1}$, the `transformation group
of~$F$' on the $2^{n-2}(n-1)!$ formally distinct expressions for~$F$.
First, note that the subgroup of~$G_{\rm asym}(n)$ isomorphic to~$S_n$
contains a subgroup $H\cong S_{n-3}$ that permutes the singular points
$a_1,\dots,a_{n-3}$.  Clearly $G_{\rm asym}(n)>G_{{\rm asym},0}(n)>H$.  The
solutions of~(\ref{eq:asym}) associated to the elements of~$H$ are simply
$F$~itself, with its $n-3$ argument pairs $\{(a_i,\rho_i)\}$ arbitrarily
permuted.  So $H$~provides the trivial part of the action of $G_{{\rm
asym},0}(n)$ on the family of equivalent expressions for~$F$.
When~$n\ge5$, $S_{n-3}$ is not normal in~$S_n$, so $H$~is not normal
in~$G_{{\rm asym},0}(n)$.  There is accordingly no natural way of
quotienting~out these trivial actions.

The following automorphisms act less trivially.  Consider the signed
permutations $[1_-][\infty_-]$ and~$[1_+\infty_+]$ of the set of singular
points $0,1,\infty,a_1,\dots,a_{n-3}$.  For a reason explained in the next
section, these will be called generalized Euler and Pfaff transformations,
respectively.  The former performs an index transformation at~$x=1$, and
the latter ``interchanges $x=1$ and $x=\infty$''; or more accurately, maps
$(0,1,\infty,a_1,\dots,a_{n-3})$ to
$(0,\infty,1,\frac{a_1}{a_1-1},\dots,\frac{a_{n-3}}{a_{n-3}-1})$, since it
performs the M\"obius transformation $x\mapsto\frac{x}{x-1}$.  Each of
these implicitly contains the positive $1$-cycle~$[0_+]$, i.e., each fixes
$x=0$ and performs no~index transformation there.  So each is an element
of~$G_{{\rm asym},0}(n)$, and yields an alternative expression for~$F$.
These are
\begin{align}
\nonumber
&F(0,\rho^{(0)};\,1,\rho^{(1)};\,\infty,\rho^{(\infty)}\!,{\hat\rho}^{(\infty)};\,a_1,\rho_1;\dots;a_{n-3},\rho_{n-3};\,Q;\,x)\\
\label{eq:genEuler}
&\qquad= (1-x)^{-\rho^{(1)}}F(0,\rho^{(0)};\,1,-\rho^{(1)};\,\infty,\rho^{(1)}\!+{\hat\rho}^{(\infty)}\!,\rho^{(1)}\!+\rho^{(\infty)}\!;\\
\nonumber
&\hskip2.0in a_1,\rho_1;\dots;a_{n-3},\rho_{n-3};\,Q';\,x)\\
\label{eq:genPfaff}
&\qquad= (1-x)^{-\rho^{(\infty)}} F(0,\rho^{(0)};\,1,{\hat\rho}^{(\infty)}\!-\rho^{(\infty)};\,\infty,\rho^{(\infty)}\!,\rho^{(1)}\!+\rho^{(\infty)};\\
\nonumber
&\hskip2.0in \tfrac{a_1}{a_1-1},\rho_1;\dots;\tfrac{a_{n-3}}{a_{n-3}-1},\rho_{n-3};\,Q'';\,\tfrac{x}{x-1}),
\end{align}
in which the expressions for $Q',Q''$, the transformed vectors of accessory
parameters, are not given explicitly since they are complicated; especially
the latter.

Each of $[1_-][\infty_-]$ and~$[1_+\infty_+]$ is an involution, as are the
transformations (\ref{eq:genEuler}) and~(\ref{eq:genPfaff}) of~$F$.  They
commute, as do (\ref{eq:genEuler}) and~(\ref{eq:genPfaff}).  So they
generate a group isomorphic to the Klein $4$-group
$\mathbb{Z}_2\times\mathbb{Z}_2$.  Besides the identity element, this group
contains the signed permutation $[1_-\infty_-]$.  By~examination, the
corresponding transformation is the same as~(\ref{eq:genPfaff}), with
$\rho^{(\infty)}$ and ${\hat\rho}^{(\infty)}$ interchanged on both sides,
and also with the sixth and seventh arguments of the right-hand~$F$ (the
exponents at~$\infty$) interchanged.  This transformation is formally
different from~(\ref{eq:genPfaff}) but is equivalent, since the order of
the two exponents at~$\infty$ is not significant on either side.

This example suggests that for computational purposes, it may be useful to
extend the automorphism group to~$\mathcal{B}_n$, and the transformation
group of~$F$ to~$\mathcal{B}_{n-1}$, by including interchange of the two
exponents at infinity.  In annotated cycle notation this interchange is the
involution~$[\infty_-]$, which acts on~$\mathcal{D}_n$ by conjugation.  So
$\mathcal{B}_n\cong \mathcal{D}_n\!\rtimes\mathbb{Z}_2$.  To~see how
computation in~$\mathcal{B}_{n-1}$, rather than~$\mathcal{D}_{n-1}$, can
facilitate the generation of the full set of $2^{n-2}(n-1)!$
transformations of~$F$, consider
\begin{equation}
\begin{split}
[1_-][\infty_-] &= ([\infty_-])^{-1}[1_+\infty_+]([\infty_-])\cdot[1_+\infty_+]\\
&= [1_-\infty_-]\cdot [1_+\infty_+],
\end{split}
\end{equation}
as an equality between transformations of~$F$.  The left side is the
generalized Euler transformation~(\ref{eq:genEuler}), and the right side is
the composition of two transformations: the generalized Pfaff
transformation~(\ref{eq:genPfaff}), which acts first, and the conjugated
Pfaff transformation mentioned above.  So although the generalized Pfaff
transformation is an involution, by computing its {\em twisted\/}
composition with itself, one obtains the generalized Euler transformation.
The two transformations are subtly related.

\section{The Hypergeometric and Heun Cases}

The hypergeometric and Heun equations, (\ref{eq:hyperg})
and~(\ref{eq:Heun}), are the $n=3$ and $n=4$ cases of the asymmetrically
reduced equation~(\ref{eq:fuchsian1}) on~$\mathbb{P}^1(\mathbb{C})$, with
the singular points taken to include~$0,1,\infty$.  These two equations
would be written as
\begin{align}
&\mathcal{E}_{\rm asym}(0,1-c;\,1,-a-b+c;\,\infty,a,b),\\
&\mathcal{E}_{\rm asym}(0,1-\gamma;\,1,1-\delta;\,\infty,\alpha,\beta;\,a,1-\epsilon;\,q),
\end{align}
in the notation of \S\,\ref{sec:basic}\,--\,\S\,\ref{sec:key}.

For these equations the canonical local solution~$F$, i.e., the Frobenius
solution at~$x=0$ with zero exponent, is defined thus.  If
$c$~(resp.~$\gamma$) is a nonpositive integer, then $F$~will typically
display logarithmic behavior at~$x=0$; that case is not discussed further
here.  If $c$~(resp.~$\gamma$) is {\em not\/} a nonpositive integer, then
$F$~will be analytic at $x=0$ with a locally convergent series expansion
$\sum_{k=0}^\infty c_kx^k$.  In the hypergeometric case $F$~will be the
Gauss function ${}_2F_1(a,b;c;\cdot)$, defined by
\begin{equation}
(k+1)(k+c)\,c_{k+1} - (k+a)(k+b) \,c_k = 0,\qquad k\ge0,
\end{equation}
and $c_0=1$.  The series $\sum_{k=0}^\infty c_kx^k$ will converge
on~$\left|x\right|<1$ to~${}_2F_1(a,b;c;x)$.  In the Heun case $F$~will be
the local Heun function $\Hl(a,q;\alpha,\beta,\gamma,\delta;\cdot)$,
defined by~\cite{Ronveaux95}
\begin{multline}
\label{eq:heuncoeffs}
(k+1)(k+\gamma)a\, c_{k+1} - \bigl\{k
\left[\,(k+\gamma+\delta-1)a + (k+\gamma+\epsilon-1)\,\right]
+q\bigr\}\,c_k \\
+ (k+\alpha-1)(k+\beta-1)\, c_{k-1} = 0,\qquad k\ge0,
\end{multline}
and the initializations $c_{-1}=0$, $c_0=1$.  The Heun series
$\sum_{k=0}^\infty c_kx^k$ will converge on the disk
$\left|x\right|<\min(1,\left|a\right|)$
to~$\Hl(a,q;\alpha,\beta,\gamma,\delta;x)$.  The exponent
parameter~$\epsilon$ is not indicated since it is constrained to equal
$\alpha+\beta-\gamma-\delta+1$, as noted above.

The technique of Proposition~\ref{prop:generate} yields $24$~formally
distinct local solutions of the hypergeometric equation, indexed by $g\in
G_{\rm asym}(3) \cong\mathcal{D}_3$.  Each equals ${}_2F_1$~with a
projectively transformed argument, multiplied by powers of $x$ and~$x-1$.
It also yields $192$~formally distinct local solutions of the Heun
equation, indexed by~$g\in G_{\rm asym}(4) \cong\mathcal{D}_4$.  Each
equals $\Hl$~with a projectively transformed argument, multiplied by powers
of $x$,~$x-1$, and~$x-a$.  The solutions are listed in Tables
\ref{tab:hyperg} and~\ref{tab:heun} respectively, normalized according to
the convention introduced in~\S\,\ref{sec:key}.

\begin{table}
\caption{Kummer's 24 local solutions of the hypergeometric equation,
indexed by~$g\in\mathcal{D}_3$.  They are partitioned into $6$~sets of
$4$~formally distinct but equivalent expressions.}
\begin{tabular}{|l|l|}
\hline\hfil $g\in\mathcal{D}_3$ &\hfil Expression corresponding to~$g$\\
\hline\hline
\multicolumn{2}{|c|}{\em Frobenius solution\/ {\rm \#1} at\/~$x=0$ {\rm(}characteristic exponent zero\/{\rm)}}\\\hline
$[0_+][1_+][\infty_+]$ & ${}_2F_1(a,b;c;\,{x})$\\
$[0_+][1_-][\infty_-]$ & ${(1-x)}^{-a-b+c}\,{}_2F_1(c-a,c-b;c;\,{x})$\\\hline
$[0_+][1_+\infty_+]$   & ${(1-x)}^{-a}\,{}_2F_1(a,c-b;c;\,\frac{x}{x-1})$\\
$[0_+][1_-\infty_-]$   & ${(1-x)}^{-b}\,{}_2F_1(c-a,b;c;\,\frac{x}{x-1})$\\\hline\hline\multicolumn{2}{|c|}{\em Frobenius solution\/ {\rm \#2} at\/~$x=0$ {\rm(}characteristic exponent\/ $1-c$\/{\rm)}}\\\hline
$[0_-][1_+][\infty_-]$ & ${x}^{1-c}\,{}_2F_1(b-c+1,a-c+1;2-c;\,{x})$\\
\hline
\end{tabular}
\label{tab:hyperg}
\end{table}

\begin{landscape}
\begin{longtable}{|l|p{6.65in}|}
\caption{The 192 local solutions of Heun's equation, indexed
 by $g\in\mathcal{D}_4$.    They are partitioned into $8$~sets of
$24$~formally distinct but equivalent expressions.}\\
\hline
\endfirsthead
\hline
\endlastfoot
\hfil $g\in\mathcal{D}_4$ &\hfil expression corresponding to~$g$\\
\hline\hline\multicolumn{2}{|c|}{{\em Frobenius solution\/ {\rm\#1} at\/~$x=0$ {\rm(}characteristic exponent zero\/{\rm)}}}\\*\hline
$[0_+][1_+][a_+][\infty_+]$ & $\Hl({a},{q};\,\alpha,\beta,\gamma,\delta;\,{x})$\\*
$[0_+][1_-][a_+][\infty_-]$ & ${(1-x)}^{1-\delta}\Hl({a},{q-(\delta-1) \gamma a};\,\beta-\delta+1,\alpha-\delta+1,\gamma,2-\delta;\,{x})$\\*
$[0_+][1_+][a_-][\infty_-]$ & ${(1-\frac{x}{a})}^{-\alpha-\beta+\gamma+\delta}\Hl({a},{q-\gamma (\alpha+\beta-\gamma-\delta)};\,-\alpha+\gamma+\delta,-\beta+\gamma+\delta,\gamma,\delta;\,{x})$\\*
$[0_+][1_-][a_-][\infty_+]$ & ${(1-x)}^{1-\delta}{(1-\frac{x}{a})}^{-\alpha-\beta+\gamma+\delta}\Hl({a},{q-\gamma [(\delta-1) a+\alpha+\beta-\gamma-\delta]};\,-\beta+\gamma+1,-\alpha+\gamma+1,\gamma,2-\delta;\,{x})$\\\hline
$[0_+][1_+\infty_+][a_+]$   & ${(1-x)}^{-\alpha}\Hl(\frac{a}{a-1},\frac{-q+\gamma \alpha a}{a-1};\,\alpha,\alpha-\delta+1,\gamma,\alpha-\beta+1;\,\frac{x}{x-1})$\\*
$[0_+][1_-\infty_-][a_+]$   & ${(1-x)}^{-\beta}\Hl(\frac{a}{a-1},\frac{-q+\gamma \beta a}{a-1};\,\beta-\delta+1,\beta,\gamma,-\alpha+\beta+1;\,\frac{x}{x-1})$\\
\end{longtable}
\end{landscape}

The tabulated solutions are partitioned respectively into $6$~sets of~$4$
and $8$~sets of~$24$.  These are sets of equivalent expressions for the~$6$
(resp.~$8$) Frobenius solutions, each defined in a neighborhood of one of
the $3$ (resp.~$4$) singular points.  They are further partitioned by the
M\"obius transformation performed on the independent variable~$x$, which is
in one-to-one correspondence with the {\em unannotated\/} cycle
representation of~$g$.  The M\"obius transformations appearing in Tables
\ref{tab:hyperg}~and~\ref{tab:heun} belong to the Kummer
group~(\ref{eq:kummergroup}) and Heun group~(\ref{eq:heungroup}),
respectively.  Although the Heun group has order~$24$, only $6$~transformed
values~$a'$ of the Heun parameter~$a$ appear in Table~\ref{tab:heun},
in~all.  This agrees with the remarks made in Example~\ref{ex:heun} on the
size of generic $a$-orbits.  If $a\in\{-1,\frac12,2\}$ or
$a\in\{\frac12\pm{\rm i}\frac{\sqrt3}2\}$, the number of distinct
values~$a'$ appearing in the $192$ solutions will be smaller: $3$~and~$2$,
respectively.

The expressions in Table~\ref{tab:hyperg} are the $24$~expressions of
Kummer~\cite{Kummer1836}, which can be found in many modern works
\cite{Abramowitz64,Dwork84,Erdelyi53,Gray2000,Poole36,Prosser94}.  But the
explicit indexing by $g\in\mathcal{D}_3$ is~unusual.
(Cf.~\cite{Prosser94}.)  Nearly all the $192$ expressions in
Table~\ref{tab:heun} have never appeared in~print before, and the indexing
by $g\in\mathcal{D}_4$ is entirely new.  The most important Frobenius
solution in each of Tables \ref{tab:hyperg} and~\ref{tab:heun} is the
first, since it provides the possible transformations of~${}_2F_1$,
resp.~$\Hl$.  In the sense of Corollary~\ref{cor:transformationgroup}, the
transformation group of~${}_2F_1$ is the order-$4$ group $G_{{\rm
asym},0}(3)\cong\mathcal{D}_2$ (i.e., the Klein $4$-group
$\mathbb{Z}_2\times\mathbb{Z}_2$), and that of~$\Hl$ is the order-$24$
group $G_{{\rm asym},0}(4)\cong\mathcal{D}_3\cong S_4$.

The transformations of~${}_2F_1$ labeled by $[1_-][\infty_-]$,
$[1_+\infty_+]$ in Table~\ref{tab:hyperg} are classical transformations due
to Euler and Pfaff~\cite{Andrews99,Prosser94}.  This explains why their
counterparts, for any~$n$, were called generalized Euler and Pfaff
transformations in~\S\ref{sec:key}.  It~is well known that Euler's
transformation of~${}_2F_1$ can be obtained by iterating
Pfaff's~\cite[\S\,2.2]{Andrews99}, but the iteration includes an
interchange of the exponents~$a,b$ at infinity, which induces a twisted
composition in~$\mathcal{D}_2$, or a composition in the order-$8$ dihedral
group~$\mathcal{B}_2$, rather than an ordinary composition
in~$\mathcal{D}_2$.  This generalizes to any~$n\ge3$.  The transformations
labeled by $[1_-][\infty_-]$, $[1_+\infty_+]$ in Table~\ref{tab:heun},
namely
\begin{align}
\nonumber
&\Hl({a},{q};\,\alpha,\beta,\gamma,\delta;\,{x})\\
\label{eq:HeungenEuler}
&\qquad= {(1-x)}^{1-\delta}\Hl({a},{q-(\delta-1) \gamma a};\,\beta-\delta+1,\alpha-\delta+1,\gamma,2-\delta;\,{x})\\
\label{eq:HeungenPfaff}
&\qquad= {(1-x)}^{-\alpha}\Hl(\tfrac{a}{a-1},\tfrac{-q+\gamma \alpha a}{a-1};\,\alpha,\alpha-\delta+1,\gamma,\alpha-\beta+1;\,\tfrac{x}{x-1}),
\end{align}
are the generalizations of the Euler and Pfaff transformations to the Heun
equation.  Snow~\cite[p.~95]{Snow52} obtained the former transformation
of~$\Hl$ from the latter by what amounts to a twisted composition
in~$\mathcal{D}_3$, or a composition in~$\mathcal{B}_3$.  Transformed
values of the accessory parameter~$q$ are painfully difficult to compute
by~hand, so for hand computation, any shortcut is valuable.

The group structure of $\mathcal{D}_3\cong S_4$ supplies many such
shortcuts, as does the structure of~$\mathcal{B}_3$, the order-$48$
extended transformation group of~$\Hl$ that includes the interchange of the
exponents~$\alpha,\beta$ at infinity.  (The $24$~equivalent expressions
for~$\Hl$ grow to~$48$ if this interchange is applied.)  Useful results
on~$\mathcal{B}_n$, such as the number of its conjugacy classes and the
maximum order of its elements, are readily available~\cite{Baake84}.  For
example, the maximum order of the elements of~$\mathcal{B}_3$ is~$6$,
though the maximum order of the elements of~$\mathcal{D}_3\cong S_4$
is~$4$.  By examination, one of the elements of order~$6$
is~$[1_+a_-\infty_+]\in\mathcal{B}_3\setminus\mathcal{D}_3$.  In
\begin{equation}
[1_+a_-\infty_+] = [1_+a_-\infty_-] \cdot [\infty_-] = [\infty_-] \cdot
[1_+a_+\infty_+],
\end{equation}
the right-multiplication by~$[\infty_-]$ interchanges $\alpha,\beta$, and
the left-multiplication by $[\infty_-]$ interchanges the third and fourth
arguments of~$\Hl$, which is not the same thing.  Regardless of how
$[1_+a_-\infty_+]$ is factored, it follows from Table~\ref{tab:heun} that
the corresponding transformation of~$\Hl$ is
\begin{align}
\nonumber
&\Hl({a},{q};\,\alpha,\beta,\gamma,\delta;\,{x})\\
&\qquad= {(1-\tfrac{x}{a})}^{-\alpha}\Hl(\tfrac{1}{1-a},\tfrac{q-\gamma \alpha}{a-1};\,-\beta+\gamma+\delta,\alpha,\gamma,\alpha-\beta+1;\,\tfrac{x}{x-a}).
\end{align}
Since $[1_+a_-\infty_+]$ generates a cyclic group of order~$6$, iterating
this will generate six transformations of~$\Hl$, including the identity.
That~is, up~to interchange of~$\alpha,\beta$, it~will generate exactly
one-fourth of the $24$ transformations of~$\Hl$, i.e., one-fourth of the
$24$~equivalent expressions for the first Frobenius solution in
Table~\ref{tab:heun}.

To see the consequences of the transformation group of~$\Hl$ being
isomorphic to~$S_4$, one should index the $24$~equivalent expressions
for~it not by $g\in\mathcal{D}_3$, but by the corresponding element $g'\in
S_4$.  An~isomorphism from $\mathcal{D}_3$ to~$S_4$ is not difficult to
construct, starting with a correspondence between their Coxeter generators.
Such an isomorphism appears in Table~\ref{tab:iso}, with each $g'\in S_4$,
the symmetric group on~$1,2,3,4$, given in cycle notation.  Since $S_4$~has
no~outer automorphisms, any isomorphism from $\mathcal{D}_3$ to~$S_4$ will
be the same as this, up~to renumbering of~$1,2,3,4$.  This isomorphism has
an elegant interpretation that can be traced back to Klein.  The
group~$\mathcal{D}_3$ can be understood as acting on a cube, as its group
of rotational symmetries.  The cube's three principal axes, which pierce
the centers of opposite sides, correspond to the singular points
$1,a,\infty$.  Under the action of any $g\in\mathcal{D}_3$ these axes are
permuted, and possibly reversed, with reversal of an axis indicated by a
negative sign in the annotated cycle representation of~$g$.  Since the
action of~$\mathcal{D}_3$ on the cube is faithfully represented by its
action on the cube's four principal diagonals, which it can permute
arbitrarily, the isomorphism $\mathcal{D}_3\cong S_4$ follows.

\begin{table}
\caption{An isomorphism between $G_{{\rm
asym},0}(4)\protect\cong\mathcal{D}_3$, the transformation group of~$\Hl$,
and~$S_4$, the symmetric group on $4$~letters.}  { \hfil
\begin{tabular}{|l|p{0.5in}|}
\hline
\hfil $g\in\mathcal{D}_3$ & \hfil $g'\in S_4$ \\
\hline\hline
$[1_+][a_+][\infty_+]$ & $(\,)$     \\
$[1_-][a_+][\infty_-]$ & $(12)(34)$ \\
$[1_+][a_-][\infty_-]$ & $(13)(24)$ \\
$[1_-][a_-][\infty_+]$ & $(14)(23)$ \\\hline
$[1_+\infty_+][a_+]$   & $(12)$     \\
$[1_-\infty_-][a_+]$   & $(34)$     \\
$[1_-\infty_+][a_-]$   & $(1423)$   \\
$[1_+\infty_-][a_-]$   & $(1324)$   \\\hline
$[1_+a_+][\infty_+]$   & $(14)$     \\
$[1_+a_-][\infty_-]$   & $(1342)$   \\
$[1_-a_+][\infty_-]$   & $(1243)$   \\
$[1_-a_-][\infty_+]$   & $(23)$     \\
\hline
\end{tabular}
\hfil
\begin{tabular}{|l|p{0.5in}|}
\hline
\hfil $g\in\mathcal{D}_3$ & \hfil $g'\in S_4$ \\
\hline\hline
$[1_+a_+\infty_+]$      & $(142)$      \\
$[1_+a_-\infty_-]$      & $(134)$      \\
$[1_-a_-\infty_+]$      & $(123)$      \\
$[1_-a_+\infty_-]$      & $(243)$      \\\hline
$[1_+][a_+\infty_+]$    & $(24)$       \\
$[1_-][a_-\infty_+]$    & $(1234)$     \\
$[1_+][a_-\infty_-]$    & $(13)$       \\
$[1_-][a_+\infty_-]$    & $(1432)$     \\\hline
$[1_+\infty_+a_+]$      & $(124)$      \\
$[1_-\infty_+a_-]$      & $(234)$      \\
$[1_-\infty_-a_+]$      & $(143)$      \\
$[1_+\infty_-a_-]$      & $(132)$      \\
\hline
\end{tabular}
\hfil
}
\label{tab:iso}
\end{table}

In this representation the Pfaff and twisted Pfaff transformations,
$[1_+\infty_+]$ and $[1_-\infty_-]$, appear as $(12)$ and~$(34)$.  Each is
a $180^\circ$ rotation about an axis passing through the midpoints of two
opposite edges of the cube.  The Euler transformation $[1_-][\infty_-]$
appears as~$(12)(34)$, making it obvious that it belongs with them in a
Klein $4$-group.  It~is a $180^\circ$ rotation about the $a$-axis.  This
representation also makes it possible to visualize the extension
of~$\mathcal{D}_3$
to~$\mathcal{B}_3=\mathcal{D}_3\rtimes_\Gamma\mathbb{Z}_2$.
On~$\mathcal{D}_3$ represented as~$S_4$ in this way, the twisting action of
the involution~$[\infty_-]$, which generates the $\mathbb{Z}_2$~factor and
interchanges the exponents at infinity, turns~out to be conjugation by the
permutation~$(14)(23)$.

The group~$S_4$ is doubly generated.  For example, it is generated by the
involution~$(14)$ and the $4$-cycle $(1423)$, which by Table~\ref{tab:iso}
correspond to $[1_+a_+]$ and~$[1_-\infty_+][a_-]$ in~$\mathcal{D}_3$.
It~follows that the full set of $24$~transformations of~$\Hl$ is generated
by {\em only two\/} well-chosen transformations, such as (by
Table~\ref{tab:heun})
\begin{align}
\nonumber
&\Hl({a},{q};\,\alpha,\beta,\gamma,\delta;\,{x})\\
&\qquad= \Hl(\tfrac{1}{a},\tfrac{q}{a};\,\alpha,\beta,\gamma,\alpha+\beta-\gamma-\delta+1;\,\tfrac{x}{a})\\
&\qquad= {(1-x)}^{\beta-\gamma-\delta}{(1-\tfrac{x}{a})}^{-\alpha-\beta+\gamma+\delta}\Hl(\tfrac{a}{a-1},\tfrac{-q-\gamma [(\beta-\gamma-\delta) a-\alpha-\beta+\gamma+\delta]}{a-1};\\
\nonumber
&\hskip1.8in -\beta+\gamma+1,-\beta+\gamma+\delta,\gamma,\alpha-\beta+1;\,\tfrac{x}{x-1}).
\end{align}
This is certainly not obvious, and would be difficult to prove without an
understanding of the group structure of the set of transformations.

The preceding examples focused on the group-theoretic aspects of the
$24$~equivalent expressions for~$\Hl$, the most important of the eight
Frobenius solutions of the Heun equation.  But the $24$~expressions given
in Table~\ref{tab:heun} for each of the remaining seven, which are
bijective with the cosets of~$\mathcal{D}_3$ in~$\mathcal{D}_4$, may also
prove useful in applications.  This is because the two-point connection
problem for the Heun equation has never been solved in full generality,
unlike the corresponding problem for the hypergeometric
equation~\cite{Schafke80a}.  No~general formula expressing the pair of
Frobenius solutions defined near a nonzero singular point (i.e.,
$x=1,a,\infty$) as a combination of the pair defined near~$x=0$, one of
which is~$\Hl$, has yet been found.

\subsection*{Some notes on Table~\ref{tab:heun}} 

The $192$ solutions were obtained with {\sc Maxima}, the free version of
the {\sc Macsyma} computer algebra system.  The implementation of the
algorithm outlined in Proposition~\ref{prop:generate} required
approximately $300$ lines of code.  

The code functioned as follows.  Each $g\in G_{\rm
asym}(4)\cong\mathcal{D}_4$ was generated as a pair $(h,\sigma)\in
(\mathbb{Z}_2)^3\!\rtimes S_4$.  Here $\sigma\in S_4$~determines the change
of independent variable $x'=P_\sigma x$, where $P_\sigma$~is one of
the~$24$ M\"obius transformations in the Heun group, listed
in~(\ref{eq:heungroup}); and $h$~specifies the subsequent index
tranformation, i.e., the negatings, if~any, of exponent differences at the
singular points $x'=0,1,a'$.  The change of independent variable and the
index transformation were applied in that order to the Heun
equation~(\ref{eq:Heun}), and the transformed parameter vector
$(a',q';\alpha',\beta',\gamma',\delta')$ was programmatically extracted
from the transformed equation.  The transformed accessory parameter~$q'$
was then re-expressed in as compact a form as possible.  The simplification
of~$q'$ sometimes required factoring quadratic polynomials in the exponent
parameters.  The first instance of this in the table occurs in the solution
indexed by~$[0_-][1_-\infty_+][a_+]$.

Each solution was then output in a linear (string) format, and converted to
\LaTeX\ format for inclusion in the table.  Conversion and table
construction were largely programmatic rather than manual.  The standard
tools {\tt sed} and~{\tt awk} distributed with the GNU/Linux operating
system were employed.

{\em The chance of there being an error in Table\/~{\rm\ref{tab:heun}} is
very small.}  After the table was prepared in \LaTeX\ format, each of the
$192$ expressions in~it was programmatically converted to the syntax of the
C~programming language, and a straightforward C~program (a~`driver') that
evaluates the $24$~expressions for any of the $8$~distinct Frobenius
solutions was written.  In~this program $\Hl$~was computed as the sum of a
Heun series, the coefficients in which were calculated from the
recurrence~(\ref{eq:heuncoeffs}).  The driver was executed with randomly
chosen values for the Heun parameters $a,q,\alpha,\beta,\gamma,\delta$, and
also for~$x$ (always taken sufficiently near the appropriate singular
point, of~course).  For any choice of these seven parameters, the
$24$~values of the specified Frobenius solution produced by the driver
always agreed with one another.  This provides strong evidence of the
correctness of Table~\ref{tab:heun}, {\em in~toto\/}.

The $192$ expressions of the table were compared with the $48$ expressions
published by Heun in his original paper~\cite{Heun1889}, and those later
published by Snow~\cite[chap.~VII]{Snow52}.  Heun's table includes every
fourth line of Table~\ref{tab:heun}, in~order (the first, fifth, ninth,
etc.).  Sadly, almost all his solutions were found to be incorrect.  It~is
not clear why he made so many errors, since it is only the transformed
accessory parameter~$q'$, which he omitted from each solution, that poses
severe difficulties in hand computation.  It is quite feasible to compute
the exponent parameters of each~$\Hl$, and the prefactor if~any, by hand,
by manipulating Riemann P-symbols.  He had especial problems with
prefactors.  For example, each of the final $12$~solutions in his
$48$-solution table has an erroneous prefactor (either $x^\alpha$
or~$x^\beta$), despite being otherwise correct.

Snow used an individualistic notation in which the exponent parameters
$\delta,\epsilon$ are interchanged, and the accessory parameter is the
negated quantity $b=-q$.  To facilitate comparison with his results, a
version of Table~\ref{tab:heun} using his notation was prepared.  (This
alternative version is available on~request from the author.)  His
$25$~solutions, each of which includes the value of~$b'$, proved to be far
more reliable than Heun's.  No~errors were found, and only a single
misprint.  In~his formula VII(21c), which gives the solution listed
as~$[\infty_+0_+1_+][a_+]$ in Table~\ref{tab:heun}, the argument of his
local Heun function should be $1/(1-z)$ rather than~$1/(z-1)$.

The several solutions supplied by the late F.~Arscott in his valuable
review of the Heun equation~\cite[Part~A]{Ronveaux95} were also examined.
In his (2.2.5\,I)--(2.2.5\,VIII) he gives the parameter vectors
$(a',q';\alpha',\beta',\gamma',\delta')$ of the eight solutions indexed by
even-signed permutations of the form $[0_\pm][1_\pm][a_\pm][\infty_\pm]$.
His transformed parameters are almost entirely correct, but there is one
misprint.  The term~`$1$' in the value for~$q'$ in~(2.2.5\,VI), which
corresponds to the solution~$[0_-][1_+][a_-][\infty_+]$ in
Table~\ref{tab:heun}, should be~`$2$'.  This misprint has already been
noted~\cite{Ronveaux2003}.  Also, sign errors were found in
(3.4.9)--(3.4.10).  Due~to these, his formula (3.4.11), which is the
Pfaff-transformed~$\Hl$ indexed in Table~\ref{tab:heun}
by~$[0_+][1_+\infty_+][a_+]$, is unfortunately incorrect.  It~is missing
its prefactor as~well, but that surely counts as a misprint.

A few additional local solutions of the Heun equation, expressed in~terms
of~${\it Hl}$, can be found scattered in the literature.  For instance,
Wakerling~\cite[\S\,IV]{Wakerling49} derived the solution indexed
by~$[\infty_+0_+][1_+][a_+]$ by manually applying to the Heun equation the
change of independent variable $x'= 1/x$, etc.  Her expression for the
transformed accessory parameter is correct.  Schmitz and Fleck, who were
the first in recent times to call attention to the problems with Heun's
original list, carefully worked~out the solutions indexed by
$[\infty_+0_+1_+a_+]$ and~$[\infty_-0_+1_+a_-]$ for their
paper~\cite[\S\,6.3.1]{Schmitz94}.  It~is painful to report that their
Eq.~(102), which gives the corresponding accessory parameters, contains
misprints (missing exponents).

In the course of writing the present paper two further lists of local
solutions of the Heun equation came to light, neither of which is cited in
standard references.  One is in the 1898 inaugural dissertation of Karl
Franz~\cite{Franz1898}, who derived $35$ of the $192$, though like Heun he
omitted accessory parameters.  He used an idiosyncratic parametrization: in
our (and Heun's) notation, his exponent parameters
$\alpha,\beta,\gamma,\gamma_1$ equal $\alpha,\beta,\gamma+\epsilon,\gamma$.
His $35$ solutions appear to be correct, so far as they go.

The other list is in the unpublished Ph.D. dissertation of Shun-Teh
Ma~\cite{Ma34}, which has languished in the Berkeley stacks for many
decades.  Though Ma evidently did not have a full group-theoretic
understanding of the transformation theory of the Heun equation, he
attempted to work~out all $192$ expressions for its Frobenius solutions,
{\em including the accessory parameters\/}, by hand computation.  This
doubtless required prodigious labor on his part.  Unfortunately a careful
comparison of his $192$ expressions with Table~\ref{tab:heun} reveals that
$64$~have incorrect exponent parameters or prefactors, and $156$~have
incorrect accessory parameters.  Of~his $192$ expressions, only $31$~have
neither sort of error.

\subsection*{A final note}

The number of local solutions of the Heun equation of the form considered
in this paper, viz., $\Hl$~with a projectively transformed argument,
multiplied by powers of $x$, $x-1$, and~$x-a$, has been variously given in
the literature.  Whittaker and Watson~\cite{Whittaker27} give the usual
figure of~$192$, which can surely be traced back to Heun.  But
Ince~\cite[p.~394]{Ince56} gives a figure of~$64$, i.e., $192/3$, and
Arscott~\cite[Part~A]{Ronveaux95} gives a figure of~$96$, i.e., $192/2$.
The variations evidently arise from differing interpretations of what it
means for two expressions for the same Frobenius solution to be distinct.
We~have been unable to find a convincing explanation for Ince's figure, but
Arscott explained the origin of his.  The involutory Heun transformation
corresponding to $[0_+][1_+a_+][\infty_+]$ in Table~\ref{tab:heun}, i.e.,
\begin{equation}
\Hl({a},{q};\,\alpha,\beta,\gamma,\delta;\,{x}) = \Hl(\tfrac{1}{a},\tfrac{q}{a};\,\alpha,\beta,\gamma,\alpha+\beta-\gamma-\delta+1;\,\tfrac{x}{a}),
\end{equation}
is unique in that its right-hand side includes no prefactor.  Applying it
to any of the $192$ expressions of the table is quite easy.  Doing so
yields another such expression, so as Arscott says, the $192$ are linked in
pairs.  In group-theoretic language, his $96$ `distinct' solutions are
indexed not by the elements of~$\mathcal{D}_4$, but rather by the cosets
in~$\mathcal{D}_4$ of the two-element group generated
by~$[0_+][1_+a_+][\infty_+]$.  Left and right cosets are different,
in~general, and he did not specify which pairing he preferred.


\end{document}